\documentclass[preprint]{elsarticle}
\usepackage{amsmath}
\usepackage{amssymb}
\usepackage{amsthm, amscd}
\newtheorem{thm}{Theorem}[section]
\newtheorem{prop}[thm]{Proposition}
\newtheorem{lem}[thm]{Lemma}
\newtheorem{cor}[thm]{Corollary}

\theoremstyle{remark}
\newtheorem{rem}[thm]{Remark}

\newcommand{\la}{\mathcal}
\newcommand{\mbb}{\mathbb}
\newcommand{\fra}{\mathfrak}
\newcommand{\g}{\fra{g}}
\newcommand{\h}{\fra{h}}

\newcommand{\bbf}{\mathbf}
\newcommand{\nrc}{\la{N}(R,r,c)}
\newcommand{\nsc}{\la{N}(S,r,c)}
\newcommand{\nn}{\la{N}}
\newcommand{\Z}{\mbb{Z}}

\newcommand{\U}{\fra{U}}

\usepackage{fancyhdr}
\begin{document}
\title{Groups elementarily equivalent to  a free nilpotent group of finite rank}
\author{Alexei G. Myasnikov\fnref{fn1}}
\ead{amiasnikov@gmail.com}
\address{Dept. Math. Stat., McGill Univ., Montreal, Quebec, Canada}
\author{Mahmood Sohrabi \fnref{fn2}}
\ead{m\_sohrabi@yahoo.com}
\address{Dept. Math. Stat., McGill Univ.,
Montreal, Quebec, Canada}

\fntext[fn1]{Partially supported by Canada Research Chair grant No. 950-201330}
\fntext[fn2]{Partially supported by NSERC grant No. RGP 298965-2005}

\begin{abstract}
In this paper we give a complete algebraic description of groups elementarily equivalent to a given free nilpotent group of finite rank.

\end{abstract}

\maketitle

\section{\protect Introduction}

\subsection{Elementary classification problem in groups}

Elementary (first-order) classification of algebraic structures goes back to the works of A. Tarski and A. Malcev. In general, the task is to characterize, in somewhat algebraic terms, all algebraic structures  elementarily equivalent to a given one. Recall, that two algebraic structures $\mathcal{A}$ and $\mathcal{B}$ in a language $L$ are elementarily equivalent ($\mathcal{A} \equiv \mathcal{B}$) if they satisfy precisely the same first-order sentences in $L$.

The first remarkable results on elementary classification of groups is due to W.Szmielew - she classified elementary theories of abelian groups in terms of ``Szmielew'' invariants \cite{Szmielew} (see also \cite{eklof,Monk, Baur}).  For non-abelian groups, the main inspiration, perhaps,  was  the long-standing  Tarski problem whether free non-abelian groups of finite rank are elementarily equivalent or not. It was  recently solved in the  affirmative in \cite{KM3,Sela6}. In contrast, free solvable (or nilpotent) groups of finite rank are elementarily equivalent if and only if they are isomorphic (A. Malcev \cite{Malcev1}).  Indeed, in these cases the abelianization $G/[G,G]$ of the group $G$ (hence the rank of $G$) is definable  (interpretable)   in $G$  by first-order formulas, hence the result.

In  \cite{malcev3} A. Malcev described elementary equivalent  classical linear groups. He showed that if $\mathcal{G} \in \{GL,PGL,SL,PSL\}$, $n, m \geq 3$, $K$ and $F$ are fields of characteristic zero, then $\mathcal{G}(F)_m \equiv \mathcal{G}(K)_n$ if and only if $m = n$ and $F \equiv K$. It turned out later that  this type of  results can be obtained via   ultrapowers by means of  the theory of abstract isomorphisms of such groups.   In this approach one argues that if the groups $\mathcal{G}(F)_m$ and $\mathcal{G}(K)_n$  are elementarily equivalent then their ultrapowers over a  non-principal ultrafilter $\omega$ are isomorphic. Since these  ultrapowers are again  groups of the type $\mathcal{G}(F^*)_m$ and $\mathcal{G}(K^*)_n$ (where $F^*$ and $K^*$ are the corresponding ultrapowers of the fields) the result follows from the description of abstract isomorphisms of such groups (which are semi-algebraic, so they preserve the algebraic scheme and the field). Similar  results    hold for many algebraic and linear groups (we refer here to a  series of papers by E. Bunina and A. Mikhalev   \cite{BM04,BM05}). On the other hand, many ``geometric'' properties of algebraic groups are just first-order definable invariants of these groups, viewed as abstract groups (no  geometry, only multiplication). For example, the geometry of a simple algebraic group is entirely determined by its group multiplication (see \cite{Z84,Poi88,Poibook}), which readily implies the celebrated Borel-Tits theorem on abstract isomorphisms of simple algebraic groups.

\subsection{On elementary classification of nilpotent groups}

In his pioneering paper \cite{malcev2} A. Malcev  showed that a ring $R$ with unit can be defined by first-order formulas in the group $UT_3(R)$ of unitriangular matrices over $R$ (viewed as an abstract group). In particular, the ring of integers $\mathbb{Z}$ is definable in the group $UT_3(\mathbb{Z})$, which is a free 2-nilpotent group of rank 2. In \cite{ershov4} Yu. Ershov  proved that the group $UT_3(\mathbb{Z})$ (hence the ring $\mathbb{Z}$) is definable in any finitely generated infinite nilpotent group $G$, which is not virtually abelian. It follows immediately that the elementary theory of $G$ is undecidable.
On the elementary classification side the main research was on M.Kargapolov's conjecture: two finitely generated nilpotent groups are elementarily equivalent if and only if they are isomorphic. In \cite{Z71}  B. Zilber gave a counterexample to the Kargapolov's conjecture. In the break-through papers \cite{MR1,MR2,MR3} A. Myasnikov and V. Remeslennikov proved that the Kargapolov's conjecture holds ``essentially'' true in the class of  nilpotent $\mathbb{Q}$-groups (i.e., divisible torsion-free nilpotent groups) finitely generated as $\mathbb{Q}$-groups. Indeed, it turned out that two such groups $G$ and $H$ are elementarily equivalent if their {\it cores} $\bar{G}$ and $\bar{H}$ are isomorphic and $G$ and $H$ either simultaneously  coincide with their cores or they do not. Here the core  of $G$ is uniquely defined as a subgroup $\bar{G} \leq G$ such that $Z(\bar{G}) \leq [\bar{G},\bar{G}]$ and $G= \bar{G} \times G_0$,  for some  abelian $\mathbb{Q}$-group  $G_0$. Developing this approach further A. Myasnikov described in \cite{alexei87,alexei90a} all groups elementarily equivalent to a given finitely generated nilpotent $K$-group $G$ over an arbitrary field of characteristic zero. Here by a $K$-group we understand P. Hall nilpotent $K$-powered groups, which are the same as $K$-points of nilpotent algebraic groups, or unipotent $K$-groups. Again, the crucial point  is that the geometric structure of the group $G$ (including the fields of definitions of the components of $G$ and their related structural constants) are first-order definable in $G$, viewed as an abstract group. Furthermore, these ideas shed some light on the Kargapolov's conjecture - it followed that two finitely generated elementarily equivalent nilpotent groups $G$ and $H$ are isomorphic, provided one of them is a {\it core} group. In this case  $G$ is a {\em core} group if $Z(G) \leq I([G,G])$, where $I([G,G])$ is the isolator of the commutant $[G,G]$. Finally, F. Oger showed in \cite{oger} that two finitely generated nilpotent groups $G$ and $H$ are elementarily equivalent if and only if they are essentially isomorphic, i.e., $G \times \mathbb{Z} \simeq H \times \mathbb{Z}$.
However, the full  classification problem for finitely generated nilpotent groups is currently wide open. In a series of papers \cite{beleg92,beleg94,beleg99} O. Belegradek completely characterized groups which are elementarily equivalent to a nilpotent group $UT_n(\mathbb{Z})$ for $n \geq 3$. It  is easy to see that (via ultrapowers) that if $\mathbb{Z} \equiv R$ for some ring $R$ then $UT_n(\mathbb{Z}) \equiv  UT_n(R)$. However, it has been shown in \cite{beleg94,beleg99} that there are groups   elementarily equivalent to $UT_3(\mathbb{Z})$ which  are not isomorphic to any group of the type $UT_3(R)$ (quasi-unitriangular groups).

\subsection{Results and the structure of the paper}

We would like now to state the main theorems proved in the paper. In the following $N_{r,c}(R)$ is a group isomorphic to the P. Hall completion of the free nilpotent group of rank $r\geq 2$ and class $c\geq 2$ over some binomial domain $R$.

\begin{thm}[Characterization Theorem]\label{charthm1111} Assume $G=N_{r,c}(R)$ and $H$ is a group. If $G\equiv H$ then $H$ is an abelian deformation of $N_{r,c}(S)$ for some ring $S$ where $R\equiv S$ as rings.\end{thm}
\begin{thm}\label{converse1} If $S$ is a ring so that $S\equiv R$ then any abelian deformation of $N_{r,c}(S)$ is elementarily equivalent to $N_{r,c}(R)$.\end{thm}
\begin{thm}\label{qnnotn1} There exists a binomial domain $R\equiv \Z$ and an abelian deformation $H$ of $N_{r,c}(R)$, for each $r\geq 2$ and $c\geq 2$, such that $H$ is not isomorphic to any Hall completion of $N_{r,c}(\Z)$.\end{thm}

In the following three subsections we introduce the notation and terminology necessary to understand the statements. The reader familiar with the theory of nilpotent groups and P. Hall completions of torsion free finitely generated nilpotent groups may skip most of~\ref{nilintro} and~\ref{hallintro} and just take note of the notations introduced. In Section~\ref{charthmsec} we prove the characterization theorem. We first give an outline of the proof in Subsection~\ref{outline} and then complete the proof in the subsequent subsection. Proofs of Theorem~\ref{converse1} and Theorem~\ref{qnnotn1} are shorter and will be included in Section~\ref{centelem} and Section~\ref{nec} respectively.

We would like to point out that when $r=2$ and $c=2$ abelian deformations of $N_{r,c}(R)$ happen to coincide with Belegradek's quasi-$UT_3$ groups (see~\cite{beleg92}). In this case all the corresponding results belong to him as mentioned earlier. In our earlier paper~\cite{MS} we gave the characterization for all free 2-nilpotent groups of finite rank. When $r>2$ or $c>2$ our main result is new to the best of our knowledge.

\subsubsection{Nilpotent groups and free nilpotent groups}\label{nilintro}We denote the lower central series of a group $G$ by
$$G=\Gamma_1(G) \geq \Gamma_{2}(G)\geq \ldots \geq \Gamma_n(G) \geq \ldots .$$
If $G$ is clear from the context we denote $\Gamma_i(G)$ by $\Gamma_i$. A group $G$ is called \textit{nilpotent} if there is a positive integer $N$ so that for all $n\geq N$, $\Gamma_n(G)=1$.  If $c$ is the least number such that
$\Gamma_{c+1}(G)=1$ then $G$ is said to be a nilpotent group of \textit{class $c$} or simply a
\textit{$c$-nilpotent} group. Let $F(n)$ be the free group on $n$ generators. A group $G$ is called a \textit{free nilpotent group of rank $n$ and class $c$} and denoted by $N_{r,c}(\Z)$, if $$G\cong F(n)/\Gamma_{c+1}(F(n)).$$

\subsubsection{P. Hall completions $N_{r,c}(R)$ of $N_{r,c}(\Z)$}\label{hallintro}
If $x,y$ is a pair of elements of a group we let $[x,y]=x^{-1}y^{-1}xy$ and call it the \textit{commutator} of $x$ and $y$. Every $G\cong N_{r,c}(\Z)$ contains an ordered tuple of elements:
$${\bf u}=(u_{11},u_{12}, \ldots , u_{1n_1},u_{21}, \ldots u_{2n_2}, \ldots, u_{c,n_c}),$$
called a \textit{Hall basic sequence} Where for each $1\leq i\leq c$,
$$\{u_{i1}\Gamma_{i+1}, u_{i2}\Gamma_{i+1}, \ldots u_{in_i}\Gamma_{i+1}\},$$
generates $\Gamma_{i}/\Gamma_{i+1}$ freely as a free abelian group. In fact each $u_{ij}$, $i>1$ is a commutator of weight $i$ in $\{u_{11}, \ldots, u_{1n_1}\}$. To avoid writing double indices let us for a moment denote the tuple above by
$${\bf u}=(u_1, u_2, \ldots , u_n).$$
Each element $g\in G$ has a unique representation:
$$g=u_1^{a_1}u_2^{a_2}\cdots u_n^{a_n} ={\bf u^a},$$
where ${\bf a}=(a_1 , \ldots ,a_n)\in \mbb{Z}^n$. Now let $h=\bf u^b$ be another element of $G$ and let $gh={\bf u^d}.$
Now if we think of $\bf a$ and $\bf b$ as tuples of $n$ integer variables then each $d_i=d_i({\bf a},{\bf b})$ is a function of $2n$ integer variables.
 On the other hand if $l$ in an integer (or on integer variable) and
$g^l=\bf u^m$
then each $m_i=m_i(l, {\bf a})$ is a function of $n+1$ integer variables $l$ and ${\bf a}$. It is a fundamental result of Philip Hall (see~\cite{hall} Section 6) that there are polynomials
$$p_i(x_1, \ldots , x_n, y_1, \ldots y_n)\in \mbb{Q}[{\bf x}, {\bf y}]$$ and
$$q_i(x_1, \ldots , x_n, y)\in \mbb{Q}[{\bf x},y]$$
called canonical polynomials associated to ${\bf u}$ such that $p_i({\bf a}, {\bf b})=d_i$ and $q_i({\bf a}, l)=m_i$. In fact the polynomials $p_i({\bf x},{\bf y})$ above are sum of integer multiples of the binomial products of the form
$$\binom{x_1}{r_1}\cdots\binom{x_i}{r_i}\binom{y_1}{s_1}\cdots \binom{y_i}{s_i}$$
and polynomials $q_i$ are integer multiples of the binomial products of the form
$$\binom{x_1}{r_1}\cdots\binom{x_i}{r_i}\binom{y}{s}$$
where the $r_i$, $s_i$ and $s$ are nonnegative integers.

Therefore if $R$ is binomial domain, i.e. $R$ is a characteristic zero integral domain
such that for all elements $r\in R$ and copies of positive integers $k=\underbrace{1+ \cdots + 1}_{k-\text{times}}$ there exists a unique solution in $R$ to the equation:
$$a (a-1)\cdots (a -k+1)= x (k!),$$
one can define a group structure (unique up to isomorphism) on $R^n$ using the polynomials $p_i$ and $q_i$. We call such a group \textit{P. Hall $R$-completion of $G$ over $R$}. In our case the P. Hall $R$-completion of $N_{r,c}(\Z)$ is denoted by $N_{r,c}(R)$. The above construction is not restricted to free nilpotent groups and can be applied to any torsion free finitely generated nilpotent group. In fact groups like $N_{r,c}(R)$ are called $R$-groups or $R$-powered groups. For details we refer the reader to~\cite{hall}.

Let us fix one more piece of notation. Whenever $\bbf{u}=(u_{11}, u_{12}, \ldots , u_{c,n_c})$ is a Hall basic sequence by $\bbf{u}_i$, $2\leq i \leq c$,
we denote the tuple $$(u_{i1}, u_{i2}, \ldots , u_{in_i}, u_{i+1,1}, \ldots , u_{c,n_c}).$$
Correspondingly by $\bbf{u}_i^{\bf a}$, where ${\bf a}=(a_{i1}, \ldots, a_{i,n_i}, a_{i+1,1}, \ldots a_{c,n_c})$ is a tuple of elements of a binomial domain $R$, we denote $u_{i1}^{a_{i1}}\cdots u_{c,n_c}^{a_{c,n_c}}$. We keep the same convention for exponents as well, i.e. if ${\bf a} = (a_{11}, a_{12}, \ldots , a_{c,n_c})$ is a tuple of elements of $R$ then by ${\bf a}_i$ we mean the sub-tuple of ${\bf a}$ with the first index greater than or equal to $i$.

\subsubsection{Abelian deformations of $N_{r,c}(R)$}
Consider the $R$-group $G=N_{r,c}(R)$. By definition there is a subset $\fra{b}=\{g_1, \ldots , g_r\}$ of $G$ with ${\bf u}$ a Hall basic sequence in $\fra{b}$ defining it as the $R$-completion of $H=N_{r,c}(\mbb{Z})$. Let $p_i$ and $q_i$ be the canonical polynomials associated to ${\bf u}$. For each $1 \leq i \leq r$ let $f^i:\oplus_{i=1}^{n_c} \times R^+\oplus_{i=1}^{n_c} R^+ \rightarrow R^+$ be a symmetric 2-cocycle. Each $f^i$ is a $n_c$-tuple of symmetric 2-cocycles $f^i_j:R^+\times R^+\rightarrow R^+$. We introduce a new product on the base set $X$ of $G$, which happens to be the set of all formal products
 $$u_{11}^{a_{11}}\cdots u_{c,n_c}^{a_{c,n_c}}={\bf u^ a},$$
 $a_{ij}\in R$. Let $g={\bf u^a}$ and $h={\bf u^b}$ be any pair of elements of this set. Now we define a product and inversion on this set as following. If $gh={\bf u^d}$ and $g^{-1}={\bf u^m}$ then
 \begin{itemize}
 \item $d_{ij}=p_{ij}({\bf a}, {\bf b})$, for all $1\leq j \leq n_i$, if $1\leq i \leq c-1$,
  \item $d_{cj} = p_{cj}({\bf a}, {\bf b})+\sum_{k=1}^r f^k_j(a_{1k}, b_{1k})$, for all $1\leq j \leq n_c$
  \item $m_{ij}=q_{ij}({\bf a}, -1)$, for all $1\leq j \leq n_i$, if $1\leq i \leq c-1$,
  \item $m_{cj}=q_{cj}({\bf a}, -1)-\sum_{k=1}^rf^k_j(a_{1k}, -a_{1k})$, for all $1\leq j \leq n_c$.
\end{itemize}
 Denote $X$ together with the operations $\cdot$ and $^{-1}$ defined above by $$N_{r,c}(R,f^1 ,f^2, \ldots ,f^r)$$ or $N_{r,c}(R,\bar{f})$. We call any group isomorphic to such a group an \textit{abelian deformation of $N_{r,c}(R)$} or a \textit{ $QN_{r,c}$-group}.

Let us recall that for an abelian group $A$ and  a group $B$ a function $f:B\times B\rightarrow A$ satisfying
\begin{itemize}
 \item $f(xy,z)f(x,y)=f(x,yz)f(y,z)$, \quad $\forall x,y,z\in B,$
 \item $f(1,x)=f(x,1)=1$, $\forall x\in B$.
\end{itemize}
is called a \textit{2-cocycle}\index{ $2$-cocycle}. If $B$ is abelian  a 2-cocycle $f:B\times B\rightarrow A$ is \textit{symmetric} if it also satisfies the identity:
$$f(x,y)=f(y,x)\quad \forall x,y \in B.$$

In Subsection~\ref{new} we shall review the correspondence between central extensions of groups and the second cohomology group and prove that abelian deformations are actually groups. We state this as a proposition below.

 \begin{prop}\label{qnisg}$N_{r,c}(R,\bar{f})$ is a group for any choice of $R$ and $\bar{f}$.\end{prop}

\begin{rem}\label{127}Let us briefly discuss the main difference between $N_{r,c}(R)$ and a $QN_{r,c}$-group over $R$. Assume $G=N_{r,c}(R)$, $\fra{b}=\{g_1, \ldots , g_r\}$ is free generating set for $G$ as an $R$-group. It is presumably known that $C_G(g_j)=g_j^R \oplus Z(G)$ for any $g_j\in \fra{b}$, where $g_j^R=\{x\in G:\exists a\in R, x=g_j^a\}$, $C_G(g)=\{x\in G:[x,g]=1\}$ and $Z(G)=\{x\in G:[x,y]=1, \forall y\in G\}$. Now consider $N_{r,c}(R,\bar{f})$ for some choice of $\bar{f}$. Then it is not hard to see that for $g_{j}\in N_{r,c}(R,\bar{f})=H$, $C_H(g_j)$ is an abelian extension (not necessarily split) of $Z(H)\cong \oplus_{i=1}^{n_c}R^+$ by $C_{H}(g_j)/Z(H)\cong R^+$ via the symmetric 2-cocycle $f^j$.\end{rem}

\section{Proof of the characterization theorem}\label{charthmsec}
In this section we give a proof of Theorem~\ref{charthm1111}. In Subsection~\ref{outline} we outline the proof. The proof will be completed in Subsection~\ref{new}.
\subsection{An outline of the proof}\label{outline}
Our approach in proving the characterization theorem resembles, in some aspects, the O.V. Belegradek's approach (see~\cite{beleg99}) in proving his characterization of groups elementarily equivalent to a unitriangular group. He finds an interpretation of the ring $R$ in the group $UT_n(R)$ in a way that any group $H\equiv UT_n(R)$ interprets a ring $S\equiv R$. Indeed he generalizes the construction due to A. Mal'cev in the case of $UT_3(R)$ to arbitrary $n$. Then he goes on to prove that the group $H$ has a structure very close to a $UT_n(S)$ except that the centralizers of standard basis elements (which happen to be abelian subgroups of $H$) are not necessarily split extensions of $S^+$ by  $Z(H)$.

 Here rather than trying to generalize Mal'cev's construction from the case of $UT_3(R)\cong N_{2,2}(R)$ to our case we take a different more global approach. Let us first state a result due to A. Myasnikov~\cite{alexei86}.
\begin{thm}[\cite{alexei86}]\label{ringinter}Let $f:M\times M\rightarrow N$ be non-degenerate full bilinear mapping of finite type. Assume $$\mathfrak{U}_R(f)=\langle R,M,N,\delta,s_M,s_N\rangle,$$ where
the predicate $\delta$ describes $f$ and $s_M$ and $s_N$ describe the actions of $R$ on the modules $M$
and $N$ respectively, and $$\mathfrak{U}(f)=\langle M,N, \delta\rangle.$$
Then there is the largest ring $P(f)$ with respect to which $f$ remains bilinear and the structure $\mathfrak{U}_{P(f)}(f)$
 is absolutely interpretable in $\mathfrak{U}(f)$. Moreover the formulas involved in the interpretation depend only on the type of $f$.\end{thm}
 Recall that an $R$-bilinear mapping $f:M\times
M\rightarrow N$ is called \textit{non-degenerate}\index{bilinear mapping!non-degenerate } in both variables if $f(x,M)=0$ or $f(M,x)=0$ implies $x=0$.
We call the bilinear map $f$, a \textit{full bilinear mapping} if $N$ is generated by $f(x,y)$, $x,y\in M$.
Let $f:M\times M\rightarrow N$ be a non-degenerate full $R$-bilinear mapping for some commutative ring $R$.
The mapping $f$ is said to have \textit{finite width}\index{bilinear mapping!of finite width} if there is a natural number $s$ such that for every $u\in
N$ there are $x_i$ and $y_i$ in $M$ we have
$$u=\sum_{i=1}^sf(x_i,y_i).$$
The least such number, $w(f)$\index{ $w(f)$}, is the \textit{width} of $f$.

A set $E=\{e_1,\ldots e_n\}$ is a \textit{complete system}\index{complete system} for a non-degenerate mapping $f$ if $f(x,E)=f(E,x)=0$ implies $x=0$. The
cardinality of a minimal complete system for $f$ is denoted by $c(f)$\index{ $c(f)$}.

\textit{Type}\index{bilinear mapping!type of}  of a bilinear mapping $f$, denoted by $\tau(f)$ \index{ $\tau(f)$}, is the pair $(w(f),c(f))$. The mapping $f$ is
said to be of finite type if $c(f)$ and $w(f)$ are both finite numbers. If $f,g:M\times M \rightarrow N$ are bilinear maps of finite type we say that the type of $g$ is less than the type of $f$ and write $\tau(g)\leq \tau(f)$ if $w(g)\leq w(f)$ and $c(g)\leq c(f)$.
The ring $P(f)$ has the following description. Let $End(M)$ denote the endomorphism ring of an abelian group $M$. We shall identify $P(f)$\index{ $P(f)$} with the subring $S\leq End(M)\times End(N)$ of all pairs $A=(\phi_1,\phi_0)$ such that for all $x,y\in M$,
 \begin{equation}\label{myrem1}f(\phi_1(x), y)=f(x,\phi_1(y))=\phi_0(f(x,y)).\end{equation}

 To make use of Theorem~\ref{ringinter} in our context we first consider the Lazard Lie ring $Lie(G)$ of the group $G$. The base abelian group of $Lie(G)$ is the direct sum $\oplus_{i=1}^c\Gamma_{i}/\Gamma_{i+1}$ and the Lie bracket on $Lie(G)$ is defined using the commutator on the group $G$. Let us introduce some notation from Lie algebra theory before proceeding.

 Let  $\g$ be any Lie algebra with respect to the bracket $[\quad, \quad]$. Define  $$\fra{g}^1=\fra{g}, \quad \fra{g}^{i+1}=[\fra{g}^i, \fra{g}]\quad \forall i\geq 1,$$
Where $[\g^i,\g]$ denotes the ideal generated by all the elements of the form $[x,y]$, $x\in \g^i$ and $y\in \g$. By $Z(\g)$ we denote the center of $\g$.
Now a \textit{free nilpotent $R$-Lie algebra of class $c$ and rank $r$}, $\nrc$, is any $R$-Lie algebra satisfying
$$\nrc\cong A(R,r)/(A(R,r))^{c+1},$$
where $A(R,r)$ is the free $R$-Lie algebra of rank $r$.

 If $\g$ is a Lie algebra the map
 $$f_{\g}:\g/Z(\g) \times \g/Z(\g) \rightarrow \g^2, \quad (x+ Z(\g) , y+Z(\g))\mapsto [x,y],$$
 is a full non-degenerate bilinear map. In case that $\g=\nn=\nrc$ one can verify that $f_{\nn}$ is of finite type. Let us point out that $\nn^c=Z(\nn)$ and that $Lie(G)\cong\nrc$ as Lie rings when $G\cong N_{r,c}(R)$. So if we can prove that $R\cong P(f_{\nn})$ and that $\nn \cong Lie(G)$ is interpretable in $G$ then we have proved that $R$ is interpretable in $G$. In fact we shall see that we can prove much more. The first main result in this direction is:

  \begin{thm}\label{R&Pf}Any element $(\phi_1, \phi_0)\in P(f_{\la{N}})$ acts by a unique element $\alpha_\phi$ of $R$ on $\la{N}/\la{N}^c$ and $\nn^2$. Moreover this correspondence is an isomorphism of rings.\end{thm}
Next we prove that
\begin{lem}\label{deflow}Each term $\Gamma_i$ of the lower central series of a finitely generated nilpotent $R$-group $G$ is absolutely definable in $G$. Moreover the same formulas define the lower central terms of any group $H\equiv G$.\end{lem}
\begin{lem}\label{grintg}The Lie ring $Lie(G)$ is absolutely interpretable in $G$. Moreover if $H\equiv G$ then $Lie(H)$ is interpreted in $H$ using the same formulas that interpret $Lie(G)$ in $G$. In particular $Lie(G)\equiv Lie(H)$.\end{lem}

\begin{lem}\label{gelemN} Let $\fra{g}$ be a Lie ring so that $\nrc\equiv \fra{g}$. Then the following statements hold.
\begin{enumerate}
\item The bilinear mapping $f_{\fra{g}}$ is absolutely interpretable in $\fra{g}$ using the same formulas that interpret $f_{\la{N}}$ in $\la{N}$. In particular $f_{\la{N}}\equiv f_{\fra{g}}$.
    \item The formulas that interpret $\fra{U}_{P(f_{\fra{g}})}(f_{\g})$ in $\fra{U}(f_{\fra{g}})$ are the same as the formulas that interpret $\fra{U}_{P(f_{\la{N}})}(f_{\la{N}})$ in $\fra{U}(f_{\la{N}})$, in particular $P(f_{\g})\equiv P(f_{\nn})$.
        \end{enumerate}\end{lem}
Now we can prove that not only the ring $R$ is interpreted in $G=N_{r,c}(R)$ but also any group $H\equiv G$ interprets a ring $S\equiv R$ and quotients of the lower central series of $H$ are free $S$-modules of the same rank as the corresponding quotients of $G$.

\begin{thm} \label{ractquot}Assume $G=N_{r,c}(R)$. Then the action of $R$ on each of the quotients $\Gamma_{i}(G)/\Gamma_{i+1}(G)$ is absolutely interpretable in $G$, i.e. the modules $$\langle  R, \Gamma_{i}/\Gamma_{i+1}, \delta_i\rangle$$ are interpretable in the group $G$. Moreover if $H\equiv G$ then there exists a ring $S\equiv R$ so that for each $i$,
$$\Gamma_i(H)/\Gamma_{i+1}(H)\cong \Gamma_i(G)/\Gamma_{i+1}(G)\otimes_{\Z}S,$$
and the formulas that interpret the action of $S$ on each $\Gamma_i(H)/\Gamma_{i+1}(H)$ are the same as the formulas that interpret the action $R$ on $\Gamma_i(G)/\Gamma_{i+1}(G)$.\end{thm}

   \textit{Proof.} Let $P=P(f_{Lie(G)})$. Notice that $\nn^2\cong \oplus_{i=2}^c\Gamma_i/\Gamma_{i+1}$. By Theorem~\ref{R&Pf}, $R\cong P$ and therefore $P$ acts on each $\Gamma_i/\Gamma_{i+1}, 2\leq i \leq c$. Hence the action of $P$ on each $\Gamma_i/\Gamma_{i+1}, 2\leq i \leq c$ is absolutely interpretable since each factor $\Gamma_i/\Gamma_{i+1}$ is so. Now consider the case $i=1$. Let us set $\g=Lie(G)$. We observe that the action of $P$ on the quotient $(\g/\g_c)/(\g^2/\g_c)$ is absolutely interpretable in $\g$. moreover the natural (group) isomorphism between $(\g/\g_c)/(\g^2/\g_c)$ and $Ab(G)$ is interpretable in $G$. This implies that the induced action of $P$ on $Ab(G)$ via the above isomorphism is absolutely interpretable in $G$. So far we have proved that all the modules $\langle P, \Gamma_i/\Gamma_{i+1}, \sigma_i \rangle$, where the predicate $\sigma_i$ describes the action of $P$ on $\Gamma_i/\Gamma_{i+1}$, are absolutely interpretable in $G$. Notice that an isomorphism of structures
$$\varphi:\U_i(P)=\langle P, \Gamma_i/\Gamma_{i+1}, \sigma_i \rangle \rightarrow \langle  R, \Gamma_{i}/\Gamma_{i+1}, \delta_i\rangle=\U_i(R)$$
is a pair $(\varphi_1,\varphi_2)$ of isomorphisms $\varphi_1:P\rightarrow R$, $\varphi_2:\Gamma_i/\Gamma_{i+1}\rightarrow \Gamma_i/\Gamma_{i+1}$ so that
\begin{equation}\label{isomulti}\varphi_2(\sigma_i(a,x))=\delta_i(\varphi_1(a),\varphi_2(x)), \quad \forall a\in P, \forall x\in \Gamma_i/\Gamma_{i+1}.\end{equation}
Now for each $i$ consider the pair $(\mu, id_i)$ where $\mu$ is the isomorphism supplied by Theorem~\ref{R&Pf} and $id_i$ is the identity map on $\Gamma_i/\Gamma_{i+1}$. It is clear by the very construction of $\mu$ that~\eqref{isomulti} is satisfied for each $i$. Thus $\U_i(P)\cong \U_i(R)$.

To prove the moreover part note that for each $i$ there are elements $\{u_{i1}, \ldots, u_{i,n_i}\}$ of $G$ whose cosets generate $\Gamma_i(G)/\Gamma_{i+1}(G)$ freely as an $R$-module. By the discussion above there exists a formula $\Phi_i(x_1, \ldots , x_{n_i})$ of the language of groups so that $G\models \Phi_i(\bar{g})$
if and only if $\bar{g}$ generate the quotient freely as $R$-module. Thus $G\models \exists \bar{x}\Phi_i(\bar{x})$ for each $i$. Therefore
$$H\models \exists \bar{x}\phi_i(\bar{x}).$$
for each $i$. By Lemma~\ref{grintg} $Lie(G)\equiv Lie(H)$ and by Lemma~\ref{gelemN}, $$S\equiv P\cong R.$$ As a consequence of Lemma~\ref{deflow} the quotients $\Gamma_i(G)/\Gamma_{i+1}(G)$ and $\Gamma_i(H)/\Gamma_{i+1}(H)$ are interpreted by the same formulas in the respective groups. So the actions of $S$ and $P$ on the corresponding quotients are interpreted by the same formulas. So $H\models \exists \bar{x}\phi_i(\bar{x})$ implies that for each $i$ there are $n_i$ elements of $H$ whose quotients generate $\Gamma_i(H)/\Gamma_{i+1}(H)$ freely as an $S$-module. This finishes the proof. \qed

The next step is to recover the action of $R$ on the group $G$ (or on some $R$ invariant subgroups of $G$) using first order formulas from the actions of $R$ on the quotients of the lower central series in a way that same formulas interpret the action of $S$ on $H$ (or on the corresponding $S$ invariant subgroups of $H$). Of course as expected the action can not be completely recovered. But we are able to prove the following statement.

 \begin{lem}\label{main}Let $G=N_{r,c}(R)$. Let $\bbf{u}$ be a Hall basic sequence for $G$. Consider the cyclic modules $u_{ij}^R=\{u_{ij}^a:a\in R\}$, viewed as structures
 $$\langle R, u_{ij}^R, \delta_{ij}\rangle,$$
 where $\delta_{ij}$ is the predicate describing the action of $R$ on $u_{ij}^R$. Then all the $u_{ij}^R$ are interpretable in $(G, \bar{u})$, where $\bar{u}=(u_{11}, \ldots , u_{1r})$, except possibly the ones generated by elements of weight 1. However when $i=1$, the action of $R$ on $C_G(u_{1j})/Z(G)$ is interpretable in $(G, u_{1j})$, for all $1\leq j \leq r$.\end{lem}

Now we are able to prove the following statement about $G$.
\begin{cor}\label{maincor}Let $G=N_{r,c}(R)$ and $\fra{b}=\{ g_1, g_2, \ldots , g_r\}$ be generating set for $G$ as an $R$-group. Let $\bf u$
be a Hall basic sequence based on $\fra{b}$. Then the following statements which are all true in $G$ can be expressed using first order formulas of the language of the enriched group $(G,u_{11}, \ldots ,u_{1r})$.
\begin{enumerate}
\item For each $1\leq i \leq c$, the set $$\{u_{i1}\Gamma_{i+1}(G), \ldots , u_{in_i}\Gamma_{i+1}(G)\}$$ generates $\Gamma_{i}(G)/\Gamma_{i+1}(G)$ freely as an $R$-module.
\item For each $1\leq j \leq r$, $C_{G}(u_{1j})$ is abelian and
$$C_G(u_{1j})/Z(G)=(u_{1j}Z(G))^R.$$
\item $u_{ij}^a u_{ij}^b = u_{ij}^{(a+b)}$ for all $a,b \in R$, if $i>1$.
\item
\begin{itemize}
\item[(a)]$[u_{ij}^ a, u_{kl}^b]={\bf u}_{i+k+1}^{{\bf t}^{ijkl}(a, b)}$ for all $a, b \in R$, if $i>1$ and $k>1$,
    \item[(b)]$[x,u_{kl}^b]={\bf u}_{k+2}^{{\bf t}^{1jkl}(a, b)}$, where $x\in C_G(u_{1j})$ and $xZ(G)=(u_{1j}Z(G))^a$, for all $a, b \in R$, if $k>1$,
    \item[(c)] $[u_{ij}^a, y]={\bf u}_{i+2}^{{\bf t}^{ij1l}(a, b)}$, where $y\in C_G(u_{1l})$ and $yZ(G)=(u_{1l}Z(G))^b$, for all $a, b \in R$, if $i>1$,
        \item[(d)] $[x,y]={\bf u}_{2}^{{\bf t}^{1j1l}(a, b)}$ where $x\in C_G(u_{1j})$ and $y\in C_G(u_{1l})$ are any elements such that $xZ(G)=(u_{1j}Z(G))^a$ and $yZ(G)=(u_{1l}Z(G))^{b}$  for all $a, b \in R$.\end{itemize}
            \end{enumerate}
\end{cor}
Thus we arrive at the following corollary:
\begin{cor}\label{1c}There is a formula $Basis(x_1, \ldots , x_r)$ of the language of groups so that if $(g_1, \ldots , g_r)\in G^r$ satisfies $Basis(\bar{x})$ then there exists a basic sequence ${\bf u}$ based on this set which satisfies statements (1.)-(4.) of Lemma~\ref{maincor}.\end{cor}
\textit {Proof.} In each item in the statement of Corollary~\ref{maincor} replace each $u_{1j}$ by the free variable $x_{1j}$ and each $u_{kl}$, $k>1$ by $x_{kl}$ where $x_{kl}$ is obtained from the $x_{1j}$ the same way as $u_{kl}$ obtained from the $u_{1j}$. One needs to notice that this way only finitely many formulas are obtained. So their conjunction produces the formula $Basis(\bar{x})$.\qed

Scrutinizing proof of Lemma~\ref{main} and Corollary~\ref{maincor} one notices that all the formulas involved in the interpretations depend only on certain logical invariants of $G$. Indeed we prove that:
\begin{lem}\label{mainlem2}Let $G=N_{r,c}(R)$ and $H$ be a group such that $G\equiv H$. Then there is a ring $S$, where $S\equiv R$ as rings, and a set $\fra{c}=\{h_1, \ldots , h_r\}$ of distinct nontrivial elements of $H$ with a Hall basic sequence ${\bf v}$ in $\fra{c}$ such that the statements (1.) - (4.) of Corollary~\ref{maincor} hold in $H$ with $R$ replaced by $S$, $u_{ij}$ replaced by $v_{ij}$ and $Z(G)$ is replaced by $Z(H)$.\end{lem}
Properties (1.)-(4.) give us all the relations we need to define a $QN_{r,c}$ group. Let us elaborate on this here.
Let
\begin{equation}\label{defI}I=_{df}\{(i,j_i)\in \mbb{N}\times \mbb{N}:1\leq i \leq c, 1\leq j_i\leq n_i\},\end{equation}
where  $c$ comes from $N_{r,c}$ and $n_i$ is the number of basic commutators of weight $i$. By definition there is a subset $\fra{b}=\{g_1, \ldots , g_r\}$ of $G=N_{r,c}(R)$ and a Hall basic sequence ${\bf u}$ in $\fra{b}$ defining it as the $R$-completion of $H=N_{r,c}(\mbb{Z})$. Then there are the canonical polynomials defining product and $R$-exponentiation in $G$. Hence for each $((i,j),(k,l))\in I\times I$ there exists a polynomial $t^{ijkl}_{rs}(x,y)\in \mbb{Q}[x,y]$, where $(r,s)\in I$, such that
$$[u_{ij}^a, u_{kl}^b]=\bbf{u}_{i+k+1}^{\bbf{t}^{ijkl}(a,b)}, \quad \forall a,b \in R.$$
We are now able to give a presentation for a $QN_{r,c}$-group.
\begin{prop}[Generators and relations for a $QN_{r,c}$ group]\label{grqn} Let ${\bf u}$ be a Hall basic sequence for $N_{r,c}(R)$. Then $N_{r,c}(R, \bar{f})$ is generated by $$\la{H}=\{u_{ij}^ a :(i,j)\in I, a\in R\}$$ and defined by the relations $\la{R}$:
 \begin{enumerate}
  \item $[u_{ij}^a, u_{kl}^b]=\bbf{u}_{i+k+1}^{{\bf t}^{ijkl}(a,b)}$, $\forall a,b \in R$, where for each $(i,j),(k,l)\in I$.
  \item $u_{ij}^a u_{ij}^b=u_{ij}^{(a+b)}$, $2\leq i \leq c$, $1\leq j=j_i \leq n_i$, $\forall a, b \in R$,
  \item $u_{1j}^a u_{1j}^b=u_{1j}^{(a+b)}\bbf{u}_c^{f^j(a, b)}$, $1\leq j \leq r$, $\forall a, b \in R$. \end{enumerate}
\end{prop}
Now the characterization theorem can be easily proved.

\textit{Proof of Theorem~\ref{charthm1111}.} To prove the statement it is enough to prove that $H$ has a presentation like the one given in Proposition~\ref{grqn} for some ring $S$ and symmetric 2-cocycles $f^j:S^+\times S^+ \rightarrow \oplus_{i=1}^{n_c}S^+$. By Lemma~\ref{mainlem2} the existence of such a presentation for $H$ is clear. We just need to remark that the symmetric 2-cocycles $f^j$ are the 2-cocycles corresponding to the abelian extension of $\oplus_{i=1}^{n_c}S^+\cong Z(H)$ by $S^+\cong C_H(u_{1j})/Z(H)$.\qed

\subsection{Proofs of the statements}\label{new}

we start by giving a brief review of the correspondence between equivalence classes (under equivalence of extensions) of central extensions of an abelian group $A$ by a
group $B$ and the group called the \textit{second cohomology group}, $H^2(B,A)$. This is included only for the convenience of the reader. The discussion is followed with a proof of Proposition~\ref{qnisg}.

By an extension of $A$ by $B$ we mean a short exact sequence of groups
$$1\rightarrow A \xrightarrow{\mu} E \xrightarrow{\nu} B \rightarrow 1,$$
where $\mu$ is the inclusion map. The extension is called abelian if $E$ is abelian and it is called central if $A\leq Z(E)$. A \textit{2-coboundary}\index{ $2$-coboundary} $g:B\times B \rightarrow A$ is a function
 defined by:
 $$\psi(xy)=g(x,y)\psi(x)\psi(y), \quad \forall x,y\in B,$$
 where $\psi:B\rightarrow A$ is a function. Any 2-coboundary is a 2-cocycle. One can make the set $Z^2(B,A)$ of all 2-cocycles  and the set $B^2(B,A)$ of all 2-coboundaries into
 abelian groups in an obvious way. Clearly $B^2(B,A)$
 is a subgroup of $Z^2(B,A)$. Let us set $$H^2(B,A)= Z^2(B,A)/B^2(B,A).$$
Assume $f$ is a 2-cocycle. Define a group $E(f)$ by $E(f)=B\times A$ as sets with the multiplication
$$(b_1, a_1)(b_2, a_2)=(b_1b_2, a_1a_2f(b_1,b_2)) \quad \forall a_1,a_2\in A,\forall b_1,b_2 \in B.$$
The above operation is a group operation and the resulting extension is central. Through out the paper we shall use the well known fact that there is a bijection between the equivalence classes of central
extensions of $A$ by $B$ and elements of the group $H^2(B,A)$ given by assigning $f+B^2(B,A)$ the equivalence class of $E(f)$.

If $B$ is abelian $f\in Z^2(B,A)$  is symmetric if and only if it arises from an abelian extension of $A$ by $B$. As it
can be easily imagined there is a one to one correspondence between the equivalent classes of abelian extensions
and the quotient group $$Ext(B,A)=S^2(B,A)/(S^2(B,A)\cap B^2(B,A))\index{ $Ext(B,A)$},$$ where $S^2(B,A)$\index{ $S^2(B,A)$} denotes the group of symmetric
2-cocycles.

For details regarding the second cohomology group we refer the reader to~\cite{robin}.

 \textit{Proof of Proposition~\ref{qnisg}.} Let $G=N_{r,c}(R)$. Set $M=\sum_{i=1}^{n_{c-1}}n_i$ and $N=\sum_{i=1}^{n_{c}}n_i$. Let
 $\bbf{a},\bbf{b}\in R^N$. Let  $\bbf{a}', \bbf{b}'\in R^M$ be the tuples of the first $M$ elements of ${\bf a}$ and $\bbf{b}$ respectively. Now for each $1\leq i \leq r$ define a function, $$g: R^M \times R^M \rightarrow R^{n_c}$$ by
 $$g(\bbf{a}', \bbf{b}')=\sum_{i=1}^rf^i(a_{1i},b_{1i}).$$
 Define
 $$k:G/Z(G)\times G/Z(G)\rightarrow Z(G),$$
 by $k(\bbf{u}^{\bbf{a}}Z(G),\bbf{u}^{\bbf{b}}Z(G))={\bf u}_c^{p_c({\bf a}', {\bf b}')+g({\bf a}', {\bf b}')}$. We need to check that $g$ defined as above is a 2-cocycle. We just remark that if ${\bf u^a  u^b}={\bf u^d}$ then $d_{1i}= a_{1i}+ b_{1i}$, $i=1,\ldots, r$. From the definition of $k$ it is obvious that $k(\bbf{u^a}Z(G),1)=k(1,{\bf u^a}Z(G))=1$ for any ${\bf a}\in R^N$. Moreover
\begin{align*}
 k({\bf u^a u^b}Z(G),{\bf u^c}Z(G)) & k({ \bf u^a}Z(G),{\bf u^b}Z(G))\\
&=\bbf{u}_c^{\sum_{i=1}^rf^i(a_{1i}+ b_{1i}, c_{1i})+f^i(a_{1i}, b_{1i})}\\
&=\bbf{u}_c^{\sum_{i=1}^rf^i(a_{1i},b_{1i}+ c_{1i})+f^i(b_{1i}, c_{1i})}\\
&=k({\bf u^a}Z(G), {\bf u^b u^c}Z(G)) k({\bf u^b}Z(G) , {\bf u^c}Z(G))
\end{align*}
for any ${\bf a,b,c} \in R^N$ which proves that $k$ is a 2-cocycle. Now it is clear that $N_{r,c}(R,\bar{f})$\index{ $N_{r,c}(R,\bar{f})$} is the central extension of $Z(G)$ by $G/Z(G)$ via the 2-cocycle $k$.
\qed

Before starting to complete the proof of characterization theorem we would like to remark that all the properties of free nilpotent Lie algebras used here are consequences of the structure theory of free Lie algebras. For the later our reference is Chapter 5 of~\cite{magnus}. In particular we shall use the following statements which are direct corollaries of Theorem 5.9 and Theorem 5.10 in the reference cited above.
\begin{lem}Consider $\la{N}=\la{N}(R,r,c)$ with a free generating set $$\fra{b}=\{\zeta_1, \ldots, \zeta_r\}.$$ Then there exists a sequence
$${\bf u}=(u_{11}, u_{12}, \ldots , u_{1,n_1},u_{21}, \ldots,u_{2,n_2},\ldots , u_{c,n_c})$$
of elements of $\nn$ called the Hall basic sequence based on $\fra{b}$ generating $\nn$ as a free $R$-module, where $u_{1j}=\zeta_j$, $j=1,\ldots r$, and each $u_{ij}$, $2\leq i\leq c$, is a certain bracket of weight $i$ in $\fra{b}$. This induces a grading
$$\la{N}\cong \la{N}_1\oplus \la{N}_2 \oplus \ldots \oplus \la{N}_c,$$
on $\nn$ where each $\la{N}_i$, $i=1,\ldots , n$, is the $R$-submodule consisting of homogeneous elements of weight $i$.
\end{lem}
\begin{lem}\label{123d} Let $z\in \la{N}$ and $z=\sum_{i=1}^nz_i$ where the $z_i$ are homogeneous of pairwise different weights. If $z=0$ then $z_i=0$ for all $i=1, \ldots n$.\end{lem}
\begin{lem}\label{123c}Assume $z$ and $t$ are homogeneous elements of $\la{N}(R,r,c)$ such that $[z,t]=0$ and the sum of their weights is strictly less than $c+1$. Then $[z,t]=0$ if and only if $z$ and $t$ are linearly dependent over $R$.\end{lem}

\textit{Proof of Theorem~\ref{R&Pf}.} Notice that in this case the bilinear map $f_{\la{N}}$ has the form:

$$\begin{array}{cccc}
f_{\la{N}}:&\la{N}/\la{N}_{c} \times \la{N}/\la{N}_c & \rightarrow & \nn^2\\
 & (x+ \la{N}_c , y+\la{N}_c)& \mapsto & [x,y]
\end{array}, $$
as $\la{N}_c=Z(\la{N})$. For simplicity we drop the subscript $\la{N}$ from $f_{\la{N}}$ and refer to it as $f$. Let us denote the coset $x+\la{N}_c$ by $[x]$ for any $x\in \la{N}$.
We start by analyzing the action of $P(f)$ on $\la{N}/\la{N}_c$ in terms of elements of $R$. To do this we pick free set of generators for $\nn$ and fix a  basic sequence ${\bf u}$ in this set. Firstly we  analyze the action of $P(f)$ on $[u_{1s}]$. Let us recall that $P(f)$ is the subring of  $End(\nn/Z(\nn)) \times End(\nn^2)$ consisting of all pairs $(\phi_1,\phi_0)$ such that
\begin{equation}\label{pf}f(\phi_1(x), y)=f(x, \phi_1(y))=\phi_0(f(x,y)), \quad \forall x,y \in \nn/Z(\nn).\end{equation}So set $\phi_1([u_{1s}])=\sum_{i=1}^{c-1}\sum_{j=1}^{n_{c-1}}\alpha_{ij}[u_{ij}]$, where each $\alpha_{ij}\in R$. From~\eqref{pf} we have
$$f(\phi_1 ([u_{1s}]) , [y])=f([u_{1s}],\phi_1([y]))=\phi_0([u_{1s}, y]),$$
for all $y\in \nn$. Letting $y=u_{1s}$ we get
\begin{equation*}\begin{split}
\sum_{i=1}^{c-1}\sum_{j=1}^{n_{c-1}}\alpha_{ij}[u_{ij}, u_{1s}]&=f(\phi_1([u_{1s}]), [u_{1s}])\\
&= \phi_0([u_{1s},u_{1s}])\\
&=\phi_0(0)=0
\end{split}\end{equation*}
Now by Lemma~\ref{123d} every homogeneous component of the sum on the left hand side of the identity above has to be zero, i.e.,
\begin{equation}\label{hom}\sum_{j=1}^{n_i}\alpha_{ij}[u_{ij},u_{1s}]=0, \quad \forall i=1, \ldots ,c-1.\end{equation}
For $i=1$ all summands in $\sum_{j=1}^{n_1}\alpha_{1j}[u_{1j},u_{1s}]$, are $R$ multiples of basic elements of weight 2 except when $j=s$. This implies that $\alpha_{1j}=0$, $1\leq j \leq r$, except possibly when $j=s$. If $i\geq 2$ then~\eqref{hom} implies that $[\sum_{j=1}^{n_i}\alpha_{ij}u_{ij},u_{1s}]=0$. But since both of the elements inside the bracket are homogeneous by Lemma~\ref{123c} we have to have that $u_{1s}$ and $\sum_{j=1}^{n_i}\alpha_{ij}u_{ij}$ are linearly dependent. This is impossible unless $\sum_{j=1}^{n_i}\alpha_{ij}u_{ij}=0$ since they are homogeneous elements of different weights. This just implies that $\alpha_{ij}=0$ for all $j=1, \ldots, n_i$. Let's fix $\alpha=\alpha_{1s}$. We shall prove that $\alpha$ obtained above is the $\alpha_\phi$ in the statement of the theorem. To do this let us first show that $\phi_1([u_{1t}])=\alpha[u_{1t}]$ for any $1\leq t \leq n_1$. So assume that $t\neq s$. By the argument above there exists
 an element $\beta$ of $R$ such that $\phi_1([u_{1t}])=\beta[u_{1t}]$. So we just need to prove that $\alpha=\beta$. Now by~\eqref{pf} applied to $x=[u_{1s}]$ and $y=[u_{1t}]$ we have $\alpha[u_{1s},u_{1t}]=\beta[u_{1s},u_{1t}]$ implying the desired identity $\alpha=\beta$ since $\nn$ is a free $R$-module. Next we prove that $\phi_1$ acts by the element $\alpha\in R$ obtained above on any element $[u_{st}]$, $1\leq s\leq c-1$ and $1\leq t \leq n_{c-1}$. So assume $1< s \leq c-1$ and let $\phi_1([u_{st}])=\sum_{i,j}\alpha_{ij}[u_{ij}]$. Consider $[u_{1k}]$ and $[u_{st}]$ for $s$ and $t$ chosen above and any $1\leq k \leq
n_1$. On the one hand
$$f(\phi_1([u_{1k}]), [u_{st}]) =  [\alpha u_{1k}, u_{st}]=[u_{1k},\alpha u_{st}].$$
On the other hand
\begin{equation*}\begin{split}
f([u_{1k}], \phi_1([u_{st}]))&=f([u_{1k}], \sum_{i,j}\alpha_{ij}[u_{ij}])\\
&= [u_{1k}, \sum_{i,j}\alpha_{ij}u_{ij}]\end{split}\end{equation*}
So by \eqref{pf} and the two identities above we have
\begin{equation*}\begin{split}\sum_{\substack{i=1\\i\neq s}}^{c-1} (\sum_{j=1}^{n_i} \alpha_{ij}[u_{1k}, u_{i,j}])
 + \sum _{\substack{j=1\\j\neq t}}^{n_s} [\alpha_{sj}u_{sj}, u_{1k}]+(\alpha_{st}-\alpha)[u_{st}, u_{1k}]=0.\end{split}\end{equation*}
 Now since $\la{N}$ is a free nilpotent Lie algebra each homogeneous element in the sum above is zero. In particular
 \begin{align*}[\sum _{\substack{j=1\\ j\neq t}}^{n_s} \alpha_{sj}u_{sj}+(\alpha_{st}-\alpha)u_{st}, u_{1k}]&=[\sum _{\substack{j=1\\j\neq t}}^{n_s} \alpha_{sj}u_{sj}, u_{1k}]+[(\alpha_{st}-\alpha)u_{st}, u_{1k}]\\
 &=0.\end{align*}
Again since the elements inside the bracket on the left hand side are homogeneous one can conclude that $\alpha_{sj}=0$ if $j\neq t$ and $\alpha=\alpha_{st}$. We also get $\alpha_{ij}=0$, if $i\neq s$. Hence $\alpha_{ij}=0$ if $(i,j)\neq (s,t)$. This proves that
$$\phi_1([x])=\alpha [x],\quad \forall x\in \la{N}. $$
Now for $u_{ck}$ a basic element of weight $c$ in $\la{N}_c$ we know that $u_{ck}=[u_{is}, u_{jt}]$ for some pair $(i,j)$ such that $i+j=c$, $1\leq s\leq n_i$ and $1\leq t\leq n_j$. So by an obvious use of \eqref{pf} we can conclude that $\phi_0(u_{ck})=\alpha u_{ck}$. It is also easy to see that for $1<i< c$ and $1\leq j \leq n_i$ we have
$\phi_0(u_{ij})=\alpha u_{ij}$, i.e. $(\phi_1, \phi_0)$ acts on $\la{N}_{i}$ by $\alpha$.

Thus we have a correspondence
$$P(f)\rightarrow R, \quad (\phi_1,\phi_0)\mapsto \alpha_\phi.$$
All the properties making the correspondence an isomorphism of unital rings are easily checked by the construction of the map.
\qed

 \textit{Proof of Lemma~\ref{gelemN}.} For any Lie ring $\fra{h}$ the ideal $Z(\fra{h})$ is absolutely definable in $\fra{h}$ by the formula
 $$\Phi_{Z}(x)=\forall y [x,y]=0.$$
  Thus $\fra{h}/Z(\fra{h})$ is absolutely interpretable in $\fra{h}$. We observe that for any Lie ring $\fra{h}$, $z\in \fra{h}^2$ if and only if $z$ satisfies one the formulas:
  $$\Psi_n(x)=_{df}\exists\bar{y}, \exists \bar{z} \quad x=\sum_{i=1}^n[y_i,z_i],$$
for some $ n\in \mathbb{N}.$
  We observe that there is a positive integer $N$, where $N$ is the number of basic elements of weight $\geq 2$, such that for every positive integer $n$ one has the following:
   $$\la{N}\models \forall x (\Psi_n(x)\rightarrow \Psi_N(x)).$$
   As $\fra{g}\equiv \la{N}$ we note that the ideal $\g^2$ is defined in $\fra{g}$ by the formula $\Psi_N(x)=\Phi_2(x)$. Now to conclude the proof of (1.) we just need to notice that bilinear maps in question are defined using Lie brackets which is already in the language so the statement (1.) follows.

   To prove (2.) one observes that by (1.) $f_{\fra{g}}$ has a type less than the type of $f_{\la{N}}$.  This implies that the formulas that interpret the ring $\fra{U}_{P(f_{\fra{g}})}(f_{\g})$ in $f_{\fra{g}}$ are the same as the formulas that interpret $\fra{U}_{P(f_{\la{N}})}(f_{\la{N}})$ in $\fra{U}(f_{\la{N}})$. Hence $P(f_{\nn})\equiv P(f_{\g})$.\qed

\textit{Proof of Proposition~\ref{grqn}.} Let $H=\langle \la{H}:\la{R} \rangle$. We notice that all the relations in the statement hold in $N_{r,c}(R,\bar{f})$. So there exists a homomorphism
$$\phi:H \rightarrow N_{r,c}(R,\bar{f}),\quad u_{ij}^a \mapsto u_{ij}^a.$$
The homomorphism $\phi$ is clearly surjective. To prove injectivity we need to prove any element $x$ of $H$ can be uniquely written
in the form $x=u_{11}^{a_{11}}\cdots u_{c,n_c}^{a_{c,n_c}}={\bf u^a}$, which is called \emph{the standard form} for $x$. This is because if $1=\phi(x)={\bf u^a}$ in $N_{r,c}(R,\bar{f})$ then $a_{ij}=0$ for all $(i,j)\in I$ (see~\eqref{defI} for the definition of $I$), which implies that $x=1$. Order the set $I$ lexicographically, i.e. $(i,j)<(k,l)$ if $i<k$ or if the two conditions $i=k$ and $j<l$ hold together. Now consider the set $\la{S}$ of all final segments of $I$ and order $\la{S}$ by comparing the least elements  of its members using $<$. If $x$ is any word in  $\la{H}$ then
  $$x=u_{k_1,l_1}^{a_1}\cdots u_{k_m,j_m}^{a_m}$$ where each $(k_i,l_i)\in I$ and $a_i\in R$. Let $I_{x}$ be the final segment of $I$ whose least element is the least subscript of $u$ in $x$. If $I_x=\{(c,n_c)\}$ by multiple applications of relation (2.) $x$ can be written in the standard form. Assume any word $w$ with $I_x< I_w$ can be written in the standard form and assume that $I_x< \{(c,n_c)\}$. Let $(k,l)$ be the least element of $I_x$. By assumption $(k,l)<(c,n_c)$. So $x$ has the form
   $$x=u_{k_1,l_1}^{a_1}\cdots u_{k_i, l_i}^{a_i}u_{k,l}^{a_{i+1}}w,$$
  $0\leq i \leq m-1$, where either $w$ is the empty word or $I_x < I_w$. By hypothesis $w$ can be written in the standard form described in the induction hypothesis. If $i=0$ we are done. So assume $i > 0$. By applying either relations (2.) or (3.) finitely many times we can assume that $(k,l)< (k_i,l_i)$. Notice that if we require to use relation (3.) the word $w$ is modified to a word $w'$. But then $I_x < I_{w'}$ and we can apply the hypothesis to $w'$ to write it in the standard form. Now use
  $$u_{k_i,l_i}^{a_{i}}u_{k,l}^{a_{i+1}}=u_{k,l}^{a_{i+1}}u_{k_i,l_i}^{a_{i}}[u_{k,l}^{a_{i+1}},u_{k_i,l_i}^{a_{i}}],$$
  and relation (1.) to get
  $$x=u_{k_1,l_1}^{a_1}\cdots u_{k_{i-1}, l_{i-1}}^{a_{i-1}}u_{k,l}^{a_{i+1}}w''$$
  where $I_x<I_{w''}$. Hence the induction hypothesis can be applied to $w''$. Hence a standard inductive argument on the number of misplaced letters with respect to
the standard form yields the result.\qed

\textit{Proof of Lemma~\ref{deflow}.} Fix a generating set $X=\{g_1, \ldots , g_m\}$ for $G$ as an $R$-group. We shall use the fact that each $\Gamma_i/\Gamma_{i+1}$ is generated as an $R$-group (here as an $R$-module) by simple commutators of weight $i$ in $X$ modulo $\Gamma_{i+1}$. We proceed by a decreasing induction on $i$. Let $c$ be the nilpotency class of $G$. Then $\Gamma_c\subseteq Z(G)$. Assume that $g_{c1}, g_{c2}, \ldots ,g_{c,m_c}$ lists all the simple commutators of weight $c$ in $X$. Then any $x\in \Gamma_c$ can be written as
 $$x=\prod_{j=1}^{m_c}g_{cj}^{a_j}, \quad a_{cj} \in R $$
 However each $g_{cj}=[g_{c-1,j_k},g_{i_k}]$ where $g_{c-1,j_k}$ is some simple commutator of weight $c-1$. Since the map
$$g_{ik}^R\rightarrow \Gamma_c, \quad x\mapsto [g_{c-1,j_k}, x],$$
is an $R$-homomorphism, we have $g_{cj}^{a_j}=[g_{c-1,j_k},g_{i_k}^{a_j}]$. So any $x \in \Gamma_c$ can be written as
 $$x=\prod_{j=1}^{m_c}[g_{c-1,j_k},g_{i_k}^{a_j}].$$
 Let $$\{C_{i,j}(g_1,\ldots , g_m):1\leq i \leq c, 1\leq j \leq m_i\},$$ for some positive integer $m_i$ list all simple commutators of weight $i$ in $X$. Hence one can define $\Gamma_c$ by
 $$\Phi_c(x)=\exists \bar{y},\exists z_1,\ldots,\exists z_{m_{c-1}}(x=\prod_{j=1}^{m_{c-1}} [C_{c-1,j}(\bar{y}),z_j]).$$

 Now fix $i< c$ and assume that for all $i\leq k \leq c$ the statement is true. Now $\Gamma_i/\Gamma_{i+1}\leq Z(G/\Gamma_{i+1})$. So by a similar argument one can conclude that for any $x\in \Gamma_{i}$, there are elements $z_1, \ldots, z_{m_{i-1}}$ such that
 $$x\Gamma_{i+1}=\prod_{j=1}^{m_{i-1}} [C_{i-1,j}(\bar{g}), z_j]\Gamma_{i+1}.$$
  Set $\Phi'(x)=_{df}\exists \bar{y},\exists z_1, \ldots \exists z_{m_{i-1}}(x=\prod_{j=1}^{m_{i-1}} [C_{i-1,j}(\bar{y}), z_j])$. Therefore by induction hypothesis $\Gamma_{i}$ is defined by the following formula:
  $$\Phi_i(x)=\exists y_1,\exists y_2 (x=y_1y_2 \wedge \Phi'_i(y_1)\wedge \Phi_{i+1}(y_2)).$$

  Now assume that $H\equiv G$. Let $S^i(y_1,\ldots , y_i)=[y_1, \ldots , y_i]$. We know that $h\in \Gamma_{i}(H)$ if and only if $h$ satisfies one of the formulas:
  $$\Psi_j(x)=\exists \bar{y}^1,\exists \bar{y}^2,\ldots, \exists \bar{y}^j (x=\prod_{k=1}^jS^i(\bar{y}^k)), $$
for some $j\in \mbb{N}$.  However for every $j\in \mbb{N}$ one has
  $$G\models \forall x (\Psi_j(x)\rightarrow \Phi_i(x)).$$
  This shows that $h\in \Gamma_i(H)$ if and only if $\Phi_i(h)$. \qed

\textit{Proof of Lemma~\ref{grintg}.} Let $\Phi_i$ be the formula defining $\Gamma_i$ in $G$ obtained in the previous lemma. Set
$$A(\bar{x})=_{\textrm{df}}x_1=x_1 \wedge (\bigwedge^c_{i=2} \Phi_i(x_i)).$$
Now define the following equivalence relation $``\sim"$ on $A$:
$$\bar{x}\sim\bar{y}\Leftrightarrow \bigwedge^{c-1}_{i=1}\Phi_{i+1}(x_iy_i^{-1}) \wedge x_c=y_c.$$

Let us denote the elements of $A/\sim$ by $[\bar{x}]$. Now define the binary operations $+$ and $[\quad, \quad]$ on $A/\sim$  by
$$\Psi_1(\bar{x}, \bar{y}, \bar{z})=_{\textrm{df}}[\bar{x}]+[\bar{y}]=[\bar{z}]\Leftrightarrow  \bigwedge^{c-1}_{i=1}\Phi_{i+1}(x_iy_iz_i^{-1}) \wedge x_cy_c=z_c,$$
$$\Psi_2(\bar{x}, \bar{y}, \bar{z})=_{\textrm{df}}[[\bar{x}],[\bar{y}]]=[\bar{z}]\Leftrightarrow \bigwedge_{k=1}^{c}\Phi_{k+1}((\prod_{i+j=k}x_i^{-1}y_j^{-1}x_iy_j)z_k^{-1}).$$
Clearly the structure obtained above is $Lie(G)$. The formulas $A$, $\Psi_1$ and $\Psi_2$ provide an absolute interpretation of $Lie(G)$ in $G$.\qed

Before embarking on the proof of Lemma~\ref{main} let us state a classical theorem of P. Hall-Petresco.
If $x_1$, \ldots , $x_r$ are free generators of a free group, we define words
$$\tau_k(x_1, \ldots, x_r=\tau_k(\bar{x}),$$
inductively by:
$$x_1^n\cdots x_r^n=\tau_1(\bar{x})^n\tau_2(\bar{x})^{\binom{n}{2}}\cdots \tau_{n-1}(\bar{x})^{\binom{n}{n-1}}\tau_n(\bar{x}).$$
\begin{thm}If $F$ is the free group on generators $x_1$ , \ldots, $x_r$ then $\tau_k(\bar{x})$ belongs to $\Gamma_k(F)$. Moreover the $\tau_i$ are independent of $n$.\end{thm}
Also recall that \textit{a simple commutator} of weight 1 in a set ${x_1, \ldots ,x_n}\in G$ where $G$ is a group is any element of the set and a simple commutator of weight $n$ has the form $[g,x_i]$ for some $i=1,\ldots n$ where $g$ is simple commutator of weight $n-1$.

\textit{Proof of Lemma~\ref{main}.} We prove that the cyclic $R$-modules generated by simple commutators of weight $\geq 2$ in $\bar{u}$ are interpretable in $(G, \bar{u})$.  Since each element of the basic sequence is a fixed product of integral powers of simple commutators of the same weight the result follows. We proceed by a decreasing induction on the weight of simple commutators.

Firstly note that $R\cong P(f_G)$ by Theorem~\ref{R&Pf} since $Lie(G)\cong \nrc$ as $R$-Lie algebras. So $R$ is absolutely interpretable in $G$ since $P(f_G)$ is interpretable in $f_G$, $f_G$ is interpretable in $Lie(G)$ by Lemma~\ref{gelemN}, and finally $Lie(G)$ is interpretable in $G$ by Lemma~\ref{grintg}.  Moreover the action of $R\cong P(f_G)$ on $Z(G)=\Gamma_c$ is interpretable in $G$ by Corollary~\ref{ractquot}. Hence the cyclic modules $u_{ci}^R$ are interpretable in $G$. Fix $k$ such that $1< k < c$. Let $l$ be the dimension of the interpretation of $R$ in $G$ and $f$ be the function from the definable subset of $G$ where $R$ is defined on onto $R$. Assume the statement is true for all simple commutators of weight $i$, $k< i \leq c$. We prove the statement for elements of weight $k$. Each simple commutator of weight $k$ is of the form  $[h,g]$ where $h$ is a simple commutator of weight $k-1$ and $g$ is a basic commutator of weight 1. Pick $a\in R$ and $y\in C_G(g)$ such that $yZ(G)=g^a Z(G)$. This choice can be made by Remark~\ref{127}. Hence there exists $v\in Z(G)$ such that $y=g^a v$. Then by Hall-Petresco formula:
\begin{equation}\label{maineq}\begin{split}
[h,y]& =[h,g^a v]\\
& = [h,v][h,g^a]^v\\
&=[h,g^a]\\
&= (h^{-1}g^{-1}h)^a g^a\\
&=[h,g]^a\tau_2(h^{-1}g^{-1}h, g)^{\binom{a}{2}}\tau_3(h^{-1}g^{-1}h, g)^{\binom{a}{3}}\cdots \tau_c(h^{-1}g^{-1}h, g)^{\binom{a}{c}}\end{split}\end{equation}

Let $g'=\tau_2(h^{-1}g^{-1}h, g)^{\binom{a}{2}}\tau_3(h^{-1}g^{-1}h, g)^{\binom{a}{3}}\cdots \tau_c(h^{-1}g^{-1}h, g)^{\binom{a}{c}}$. Then $g'$ is an element of $\Gamma_{k+1}(G)$. Each $\tau_m(h^{-1}g^{-1}h,g)$ is a product of integral powers of commutators in $h^{-1}g^{-1}h$ and $g$. So there are integers $b^m_{ij}$ such that
$$\tau_m(h^{-1}g^{-1}h, g)=\prod_{i=m+k-1}^c\prod_{j=1}^{n_i}u_{ij}^{b^m_{ij}}.$$
Now the existence of the canonical polynomials associated to $\bbf{u}$ implies
the existence of polynomials $$r_{ij}(x_1, \ldots x_c, {\bf y^{k+1}}, \ldots , {\bf y^{c}})$$ where ${\bf y}^i=(y^i_{11}, \ldots, y^i_{c,n_c}),$ so that
$$g'=\prod_{(i,j)\in I}u_{ij}^{r_{ij}(\binom{a}{1}, \ldots ,\binom{a}{c}, {\bf b^{k+1}}, \ldots, {\bf b^{c}})},$$
where $r_{ij}(\binom{a}{1}, \ldots ,\binom{a}{c}, {\bf b^{k+1}}, \ldots, \bbf{b^{c}})=0$ whenever $i\leq k$.
Since actually each $r_{ij}$ is a sum of integral multiples of products of binomial coefficients there is an equation expressible in the first order language of rings so that its unique solution is $r_{ij}(\binom{a}{1}, \ldots ,\binom{a}{c}, {\bf b^{k+1}},\ldots, {\bf b^{c}})$. Now by induction hypothesis each cyclic module $u_{ij}^R$, $i> k$ is interpretable in $(G,\bar{u})$. So there exists a first order formula $\Phi (x, y_1, \ldots , y_r, z_1, \ldots , z_l)$ (note that any $u_{ij}$ is a certain commutator in $\bar{u}$) of the language of groups such that
\begin{equation*}\begin{split}
    g'& =\tau_2(h^{-1}g^{-1}h, g)^{\binom{a}{2}}\tau_3(h^{-1}g^{-1}h, g)^{\binom{a}{3}}\ldots \tau_c(h^{-1}g^{-1}h, g)^{\binom{a}{c}}\\
& \Leftrightarrow (G, \bar{u})\models \Phi (g', \bar{u}, g_1, \ldots , g_l)
\end{split}\end{equation*}
 where $f(g_1, \ldots, g_l)=a$. By Corollary~\ref{ractquot} the action of $R$ on $Ab(G)$ is interpretable in $G$ so clearly there is a formula $\Phi'$ of the language of groups so that
 $$y\Gamma_2=(g\Gamma_2)^a \Leftrightarrow G\models \Phi'(y,g,g_1, \ldots , g_l).$$
 So we have
\begin{equation}\label{simpleint}\begin{split}
 x=[h,g]^a
 &\Leftrightarrow
 (G,\bar{u})\models \exists z,y (x=[h,y]z^{-1} \wedge \Phi (z,\bar{u}, g_1, \ldots , g_l)\\
 &\quad \wedge  \Phi'(y,g,g_1, \ldots , g_l) \wedge [g,y]=1).
\end{split}\end{equation}
 Thus the formula on the right hand side of $\Leftrightarrow$ in \eqref{simpleint} interprets the action of $R$ on the abelian group $[h,g]^R$ with respect to the parameters $\bar{u}$. We notice that $g$ and $h$ chosen above are some specific commutators in $\bar{u}$.

In order to prove that the action of $R$ on $C_G(u_{1j})/Z(G))$ is interpretable in $(G, u_{1j})$ firstly we notice that $C_G(u_{1j})=u_{1j}^R\oplus Z(G)$, so the following equivalence should be clear.
$$xZ(G)= (u_{1j}Z(G))^a \Leftrightarrow x\Gamma_2(G)=(u_{1j}\Gamma_2(G))^a \wedge [x, u_{1j}]=1,$$
for all $a\in R$. But the right hand side is expressible in the first order language of the enriched group $(G,u_{1j})$. The result follows now.\qed

\textit{Proof of Corollary~\ref{maincor}.} Statement (1.) is expressible by formulas of the language of $(G, {\bf u})$ since the action of $R$ on each $\Gamma_i(G)/\Gamma_{i+1}(G)$ is absolutely interpretable in $G$. The result for (2.), (3.) and (4.) is a direct consequence of Lemma~\ref{main}. \qed

 \textit{Proof of Lemma~\ref{mainlem2}.} Let $Basis(\bar{x})$ be the formula obtained in Corollary~\ref{1c}. Since $H\equiv G$, we have
$$H\models \exists \bar{x}Basis(\bar{x}).$$
Let $(h_1, \ldots, h_r)$ be a tuple of elements of $H$ such that $H\models Basis(\bar{h})$ and ${\bf v}$ be the Hall basic sequence based on these elements. Set $S=P(f_{Lie(H)})$. As a corollary of Theorem~\ref{ractquot} statement $(1.)$ in Corollary~\ref{maincor} holds in $H$ with $u_{ij}$ replaced by $v_{ij}$.

  Moreover
 $$xZ(G)= (u_{1j}Z(G))^a \Leftrightarrow x\Gamma_2(G)=(u_{1j}\Gamma_2(G))^a \wedge [x,G]=1.$$
  So the right hand side of ``$\Leftrightarrow$'' can be used with corresponding replacements to interpret the action of $S$ on $C_H(v_{1j})/Z(H)$. This proves that Statement $(2.)$ holds in $H$ with proper replacements.

To prove that $(3.)$ and $(4.)$ are true in $H$ with proper replacements we will first prove that for $2 \leq i \leq c$ the sets
$$v_{ij}^S=\{v_{ij}^a: a \in S\}$$
are cyclic $S$-modules which are interpretable in the enriched structure $(H, \bar{v})$. To do this we observe that $v_{ij}$, $2\leq i \leq c$, are products of integral powers of simple commutators in $\{v_{1j}:1\leq j \leq r\}$ since the same relations hold between the $u_{ij}$ and the $u_{1j}$. Now using a decreasing induction on the weight of simple commutators in $\{v_{1j}:1\leq j \leq c\}$ we let Equation~\eqref{maineq} define the $S$ exponents of these simple commutators. So by the observation made above the $S$ exponents of each $v_{ij}$, $1\leq i \leq c$, $1\leq j \leq n_c$ can be defined. Now since each $u_{ij}^R$ is a cyclic module and $S$ exponentiation in $v_{ij}^S$ is defined using the action of $R$ on $u_{ij}^R$, $S$ exponentiation is actually an action and turns $v_{ij}^S$ into $S$-modules. We just remark that the $S$-module structure of each $v_{ij}^S$ is interpretable in $(H, \bar{v})$ using the same formulas that interpret the action of $R$ on $u_{ij}^R$. Moreover from the above paragraph we have that the action of $S$ on $C_H(v_{1j})/Z(H)$ is interpreted in $H$ using the same formulas that interpret the action of $R$ on $C_G(u_{1j})/Z(G)$. The final point to consider is the polynomials ${\bf t}^{ijkl}$.These polynomials make sense over any binomial domain. Since $R$ is a binomial domain and $R\equiv S$ hence is $S$. So the polynomials ${\bf t}^{ijkl}$ can be regarded to be the same if we identify the copies of $\mbb{Z}$ inside the two rings. The statement follows now.
 \qed


\section{Central extensions and elementary equivalence}\label{centelem}
The aim of this section is to prove that for any two elementarily equivalent binomial domains $R$ and $S$
$$ N_{r,c}(R,\bar{f})\equiv N_{r,c}(S,\bar{g})$$ for any symmetric 2-cocycles $f^i$ and $g^i$, $1\leq i \leq n_c$.
 \begin{lem}\label{125}The group $N_{r,c}(R)$ is absolutely interpretable in the ring $R$ and the formulas involved in the interpretation depend only on $R$ being a binomial domain.\end{lem}
 \textit{Proof.} The polynomials $p$ and $q$ provide a $\sum_{i=1}^cn_i$ dimensional interpretation of $N_{r,c}(R)$ in $R$. Indeed this object is the group of $R$ points of a nilpotent algebraic group. Since the formulas involved in the interpretation depend only on $p$ and $q$ and in turn $p$ and $q$ do not depend on $R$ as far as $R$ is a binomial domain the statement follows.\qed
 \begin{cor}If $R\equiv S$ where $R$ and $S$ are some binomial domains then $N_{r,c}(R)\equiv N_{r,c}(S)$.\end{cor}

\begin{lem}\label{isocent}
Let
$$1 \rightarrow A \rightarrow G \rightarrow B\rightarrow 1$$
be a central extension of an abelian groups $A$ by a group $B$. Let $(J,\mathcal{D})$ be an ultrafilter. Then $G^J/\mathcal{D}$ is isomorphic to a central extension of $A^J/\mathcal{D}$ by $B^J/\mathcal{D}$.\end{lem}
We omit the proof of the lemma. Its proof is completely similar to the proof of~\cite{MS}, Lemma 7.1. Otherwise the proof is elementary.

\begin{lem}\label{converse} For any choice of $\bar{f}$, $N_{r,c}(R)\equiv N_{r,c}(R,\bar{f})$.\end{lem}

\textit{Proof.} We will prove the statement using ultrapowers. We need to remark that if $R$ is a binomial domain then for any ultrafilter $(J,\la{D})$ , $R^J/\la{D}$ is also a binomial domain.

In Lemma~\ref{qnisg} we obtained a 2-cocycle $k$ so that $G=N_{r,c}(R, \bar{f})$ is a central extension of $Z(G)$ by $G/Z(G)$ via $k$.

Choose $(J, \la{D})$ so that $(R^+)^J/\la{D}$ is $\omega_1$-saturated. Now by Theorem 1.10 of~\cite{eklof} we can conclude that each $$(f^i)^{\la{D}}\in B^2((R^+)^J/\la{D}, \oplus_{i=1}^{n_c} (R^+)^J/\la{D}),$$ where $(f^i)^\la{D}$ is the obvious 2-cocycle induced by $f^i$, which in turn implies that  $\sum_{i=1}^r(f^i)^{\la{D}}\in B^2( \oplus_{i=1}^r(R^+)^J/\la{D}, \oplus_{i=1}^{n_c} (R^+)^J/\la{D})$. Now by definition of the 2-cocycle $k$, $H=N_{r,c}(R^J/\la{D})$ and $$H'=N_{r,c}(R^J/\la{D}, (f^1)^{\la{D}}, \ldots , (f^r)^{\la{D}})$$ are equivalent as extensions of $Z(H)$ by $H/Z(H)$. So in particular
$$N_{r,c}(R^J/\la{D})\cong N_{r,c}(R^J/\la{D}, (f^1)^{\la{D}}, \ldots , (f^r)^{\la{D}}).$$

On the other hand by Lemma~\ref{isocent},
$$(N_{r,c}(R))^J/\la{D}\cong N_{r,c}(R^J/\la{D}),$$ and
$$N_{r,c}(R^J/\la{D}, (f^1)^{\la{D}}, \ldots , (f^r)^{\la{D}})\cong (N_{r,c}(R, \bar{f}))^J/\la{D},$$

hence, $N_{r,c}(R)\equiv N_{r,c}(R,\bar{f})$.\qed

Thus Theorem~\ref{converse1} is a direct corollary of Lemma~\ref{125} and Lemma~\ref{converse}.

\section{A $QN_{r,c}$-group which is not $N_{r,c}$}\label{nec}
In this section we prove the existence of a $QN_{r,c}$-group over a certain ring which is not a $N_{r,c}$-group over any ring. Before that we need to have a description of abstract isomorphisms of free nilpotent Lie algebras.
\begin{lem}\label{ailem}Let $\psi:\nrc \rightarrow \nsc$ be a Lie ring isomorphism. Let $\psi_1:\fra{g}/Z(\fra{g})\rightarrow \fra{h}/Z(\fra{h})$ and $\psi_0:\g^2\rightarrow \h^2$ be the isomorphisms induced by $\psi$.Then there is a ring isomorphism $\mu:R\rightarrow S$ such that
 $$\psi_1(a(x+Z(\g))=\mu (a)\psi_1(x+Z(\g )), \quad \forall a\in R, \forall x\in \g,$$
 and
  $$\psi_0(ax)=\mu (a)\psi_0(x), \quad \forall a\in R, \forall x\in \g^2.$$
\end{lem}

 \textit{Proof.} We prove that $P(f_\g)\cong P(f_\h)$. Then Theorem~\ref{R&Pf} implies the existence of an isomorphism between $R$ and $S$. Consider the map:
$$\mu:P(f_{\fra{g}})\rightarrow P(f_{\fra{h}}), \quad (\phi_1, \phi_0)\mapsto (\psi_1\phi_1\psi_1^{-1}, \psi_0\phi_0\psi_0^{-1})).$$
Firstly we need to check if $P(f_{\h})$ is actually the target of the map defined. This is a consequence of the fact that $\psi$ is a Lie ring isomorphism. Indeed pick any $x,y\in \h$. Then,
\begin{align*}
f_{\h}(\psi_1\phi_1\psi_1^{-1}(x+Z(\h)), y+z(\h))&=\psi_0(f_{\g}(\phi_1\psi_1^{-1}(x+Z(\h)),\psi_1^{-1}(y+Z(\h)))\\
&=\psi_0\phi_0f_{\g}(\psi_1^{-1}(x+Z(\h)),\psi^{-1}(y+Z(\h)))\\
&=\psi_0\phi_0\psi_0^{-1}f_\h (x+Z(\h),y+Z(\h)).\end{align*}
The map $\mu$ being a homomorphism follows from
$$\psi_i(\phi_i+\phi'_i)\psi_i^{-1}=\psi_i\phi_i\psi_i^{-1}+ \psi_i\phi'_i\psi_i^{-1},$$
and $$\psi_i\phi_i\phi'_i\psi_i^{-1}=\psi_i\phi_i\psi_i^{-1}\psi_i\phi'_i\psi_i^{-1},$$
$i=1,2$.

One can easily check that
$$\mu':P(f_\h)\rightarrow P(f_\g), \quad (\theta_1,\theta_0)\mapsto (\psi_1^{-1}\theta_1\psi_1,\psi_0^{-1}\theta_0\psi_0),$$
is the inverse of $\mu$.

 We also denote the isomorphism obtained from $R$ to $S$ by $\mu$. Now choose $x\in \fra{g}$ and let $a \in R$. By Theorem ~\ref{R&Pf} there exists $(\phi_1,\phi_0) \in P(f_{\fra{g}})$ so that $a ( x+Z(\g)) = \phi_1(x+Z(\g))$. Then
\begin{align*}
 \psi_1(a(x+Z(\g ))&= \psi_1\phi_1(x+Z(\g))\\
 &=\psi_1\phi_1\psi_1^{-1}\psi_1(x+Z(\g))\\
 &=\mu(a) \psi_1(x+Z(\g)).
\end{align*}
One can repeat this argument for $\psi_0$ easily. This finishes the proof.\qed
\begin{lem}\label{126a}
assume that $\eta: G=N_{r,c}(R,\bar{f})\rightarrow N_{r,c}(S)=H$ is an isomorphism of groups. Then the rings $R$ and $S$ are isomorphic via a map $\mu:R\rightarrow S$. Moreover if $\eta_1:Ab(G)\rightarrow Ab(H)$ is the isomorphism induced by $\eta$ and $\eta_0:Z(G)\rightarrow Z(H)$ is the restriction of $\eta$ to $Z(G)$ then we have
$$\eta_1((x\Gamma_2(G))^a)=(\eta_1(x\Gamma_2(H)))^{\mu(a)},\quad \forall x\in G,\forall a\in R,$$
and
$$\eta(x)=(\eta(x))^{\mu(a)}, \quad \forall x\in Z(G),  \forall a\in R.$$\end{lem}

\textit{Proof.} To obtain the isomorphism $\mu:R\rightarrow S$ we may use Lemma~\ref{ailem} since $Lie(G)\cong \la{N}(R,r,c)$ and $Lie(H)\cong\la{N}(S,r,c)$, as Lie algebras.
Consider the Lie ring isomorphism $\psi: \fra{g}=Lie(G)\rightarrow Lie(H)=\fra{h}$ induced by $\eta$ and define $\psi_1$, $\psi_0$ and $\mu$ similar to the ones described in Lemma~\ref{ailem}.

  Any $x\in \g$ can be uniquely written as $(x)_1+ (x)_2$ where $(x)_1\in Ab(G)$ and $(x)_2\in \g^2$. So we obviously have that
 \begin{align*}
 \eta_1((x\Gamma_2(G))^a )&=(\psi (a(x\Gamma_2(G)))_1\\
 & =(\mu (a) \psi(x\Gamma_2(G))+z)_1, \quad \textrm{for some }z\in Z(\g)\\
 &=(\mu (a) \psi(x\Gamma_2(G)))_1\\
 &=(\eta_1(x\Gamma_2))^{\mu(a)}
\end{align*}
for all $x$ in $G$. A similar argument using $\phi_0$ and $\psi_0$ instead of $\phi_1$ and $\psi_1$ proves that
$$\eta(x^a)=\eta(x)^{\mu(a)}, \quad \textrm{for all  }x\in Z(G).$$\qed
\begin{thm}\label{ghaz}
Let $\eta:G=N_{r,c}(R,\bar{f})\rightarrow N_{r,c}(S)=H$ be an isomorphism of groups. Then for each $1\leq j \leq r$ we have that $f^j\in B^2(R^+, \oplus_{i=1}^{n_c}R^+)$, i.e. each $f^j$ is a 2-coboundary.
\end{thm}
\textit{Proof.} Let the tuple of elements ${\bf u}$ of $G$ be the one appeared in the definition of a $QN_{2,n}$ group. Set $\eta(u_{ij})=v_{ij}$ for all $(i,j)\in I$. Let $\eta_1:Ab(G)\rightarrow Ab(H)$ be the group isomorphism induced by $\eta$. By Lemma~\ref{126a} there exists an isomorphism $\mu:R \rightarrow S$ of rings so that $\eta_1(x\Gamma_2(G)^a)=(\eta_1(x\Gamma_2(G)))^{\mu(a)}$, for all $a$ in $R$ and $x$ in $G$. This implies that $\{v_{11}\Gamma_2(H), \ldots, v_{1r}\Gamma_2(H)\}$ generates $Ab(H)$ freely as an $S$-module since $\{u_{11}\Gamma_2(G), \ldots, u_{1r}\Gamma_2(G)\}$ generates $Ab(G)$ freely as an $R$-module. So $\fra{c}=\{v_{11}, \ldots , v_{1r}\}$ generates $H$ as an $S$-group. Let $
{\bf v}$ be the Hall basic sequence in $\fra{c}$ then every element $h$ of $H$ has a unique representation ${\bf v^a}$.  Set $J=\{(i,j)\in I:2\leq i \leq c\}$ and $M=\sum_{i=2}^{n_{c}}n_i$. By Lemma~\ref{126a}
$$\eta(u_{1j}^a)=v_{1j}^{\mu(a)}{\bf v}_2^{g(\mu(a))}, \quad \forall a \in R,$$
where $g=(g_{ij})_{(i,j) \in J} :S \rightarrow S^M$ is a function determined by $\eta$. Since $u_{1j}^a\in C_G(u_{1j})$ we have to have that $v_{1j}^{\mu(a)}{\bf v}_2^{g(\mu(a))}\in C_H(v_{1j})$. Remark~\ref{127} implies that $g_{ij}=0$ for all $(i,j)$ such that $2\leq i\leq n_{c-1}$. Hence one could write
$$\phi(u_{1j}^a)=v_{1j}^{\mu(a)}{\bf v}_c^{g(\mu(a))}, \quad \forall a \in R.$$
Choose two arbitrary elements $b$ and $b'$ in $S$. Then,

\begin{align*}
 v_{1j}^{b+b'}&=v_{1j}^b v_{1j}^{b'}\\
 &=\phi (u_{1j}^{\mu^{-1}(b )}){\bf v}_c^{-g(b)}\phi (u_{1j}^{\mu^{-1}(b' )}){\bf v}_c^{-g(b')}\\
 &=\phi (u_{1j}^{\mu^{-1}(b )})\phi (u_{1j}^{\mu^{-1}(b' )}){\bf v}_c^{-g(b)-g(b')}\\
 &=\phi (u_{1j}^{\mu^{-1}(b +b')}{\bf u}_c^{f^j(\mu^{-1}(b ), \mu^{-1}(b'))}){\bf v}_c^{-g(b)-g(b')}\\
 &=v_{1j}^{b+b'}{\bf v}_c^{\mu f^j (\mu^{-1}(b ), \mu^{-1}(b'))+g(b+b')-g(b)-g(b')},
 \end{align*}
where $\mu f^j=_{df}(\mu f^j_k)_{1\leq k \leq n_c}$.
The identity above clearly shows that $$\mu f^j (\mu^{-1}(-), \mu^{-1}(-))\in B^2(S^+, \oplus_{i=1}^{n_c}S^+).$$ Since $\mu$ is a ring isomorphism this implies that all $f^j$, $j=1,\ldots ,r$, are 2-coboundaries as claimed.\qed

\begin{lem}[O. V. Belegradek, \cite{beleg92}]\label{beleglem}There is a ring $R$, $R\equiv \mathbb{Z}$ such that $Ext(R^+,R^+)\neq 0$.\end{lem}
\textit{Proof of Theorem~\ref{qnnotn1}.} By Lemma~\ref{beleglem} there exists a ring $R$ such that $R\equiv \mathbb{Z}$ and $Ext(R^+,R^+)\neq0$. Then by Theorem~\ref{ghaz} there has to exist 2-cocycles $f^i:R^+\times R^+\rightarrow \oplus^{(^n_c)}_{i=1} R^+$, $1\leq i \leq n$, such that
$$H=N_{r,c}(R,f^1,\ldots , f^n) \ncong N_{r,c}(R).$$
We note that $H\ncong N_{r,c}(S)$ for any binomial domain $S$ by Theorem~\ref{ghaz}. Moreover $H \equiv G$ by Lemma~\ref{converse}.\qed

\vspace{2cm}


\end{document}